\documentclass[12pt,final]{amsart}
\usepackage{fancyhdr}
\usepackage{stmaryrd}
\usepackage{graphicx}
\usepackage{setspace}
\usepackage{verbatim}
\usepackage{scalerel}
\usepackage{amssymb}

\usepackage[a4paper,left=3cm,right=3cm,top=3cm,bottom=4.5cm]{geometry}
\usepackage[hypcap]{caption}
\usepackage{comment}
\usepackage{etaremune}
\usepackage{graphicx}
\usepackage{hyperref}
\usepackage{color}
\usepackage{pdfsync}
\usepackage{showlabels}

\theoremstyle{definition}
	\newtheorem{thm}{Theorem}[section]
        \newtheorem*{thm*}{Theorem}
        \newtheorem{Def}[thm]{Definition}
	\newtheorem*{Def*}{Definition}
	\newtheorem{rmk}[thm]{Remark}
	\newtheorem*{rmk*}{Remark}
	\newtheorem{cor}[thm]{Corollary}
	\newtheorem*{cor*}{Corollary}
	\newtheorem{prop}[thm]{Proposition}
	\newtheorem*{prop*}{Proposition}
	\newtheorem{eg}[thm]{Example}
	\newtheorem*{eg*}{Example}
	\newtheorem{lem}[thm]{Lemma}
	\newtheorem*{lem*}{Lemma}
	
        \newtheorem*{ex*}{Exercise}
        
	\newtheorem*{claim*}{Claim}
        \newtheorem{fact}[thm]{Fact}
\theoremstyle{definition}
	
	\newtheorem{soln*}{Solution}
	
	\newtheorem*{note*}{Note}
        
   \newtheorem{conv}[thm]{Convention}

\newcommand{\real}{\mathbb{R}}
\newcommand{\complex}{\mathbb{C}}

\newcommand{\ve}{\varepsilon}
\newcommand{\OO}{\mathcal{O}}
\newcommand{\Vv}{\mathcal{V}}

\newcommand{\VG}{\Gamma}
\newcommand{\Stab}{\mathrm{Stab}}
\newcommand{\res}{\mathrm{res}}
\newcommand{\val}{\mathrm{val}}

\newcommand{\ti}[1]{\tilde{#1}}
\newcommand{\ol}[1]{\overline{#1}}
\newcommand{\AAA}{\mathbb{A}}
\newcommand{\PPP}{\mathbb{P}}
\newcommand{\StackP}[2][Tag]{\cite[\href{https://stacks.math.columbia.edu/tag/#2}{#1~#2}]{stacks-project}}
\newcommand{\St}[2]{\{#1:#2\}}

\newcommand{\GL}{\mathrm{GL}}
\newcommand{\SL}{\mathrm{SL}}
\newcommand{\PSL}{\mathrm{PSL}}
\newcommand{\Gm}{\mathrm{G_m}}
\newcommand{\DefSet}[1]{\mathbf{#1}}
\newcommand{\RES}{\DefSet{RES}}
\newcommand{\VF}{\DefSet{VF}}
\newcommand{\GV}{\DefSet{G}}
\newcommand{\GO}{\GV(\OO)}
\newcommand{\GR}{\ol{\GV}}
\newcommand{\HH}{\DefSet{H}}
\newcommand{\XX}{\DefSet{X}}
\newcommand{\PP}{\DefSet{P}}
\newcommand{\YY}{\DefSet{Y}}
\newcommand{\VV}{\DefSet{V}}
\newcommand{\ZZ}{\DefSet{Z}}
\newcommand{\WW}{\DefSet{W}}
\newcommand{\BB}{\DefSet{B}}
\newcommand{\CC}{\DefSet{C}}

\newcommand{\mutube}[1]{\mu\cdot{#1}}

\newcommand{\Th}[1]{\mathrm{#1}}
\newcommand{\ACVF}{\Th{ACVF}}

\newcommand{\Lang}[1]{\mathcal{#1}}
\newcommand{\LL}{\Lang{L}}

\newcommand{\Model}[1]{\mathbb{#1}}
\newcommand{\Mm}{\Model{M}}
\newcommand{\Uu}{\Model{U}}
\newcommand{\Ll}{\Model{L}}

\newcommand{\Mat}[4]{%
  \bigl(\begin{smallmatrix}#1&#2\\#3&#4\end{smallmatrix}\bigr)}
\newcommand{\kk}{\Bbbk}
\newcommand{\ra}{\rightarrow}
\newcommand{\deft}[1]{\emph{#1}}

\title{Peterzil-Steinhorn subgroups and \(\mu\)-stabilizers in ACF}
\author{Moshe Kamensky}
\address{Department of Mathematics, Ben-Gurion University}
\email{kamenskm@math.bgu.ac.il}
\author{Sergei Starchenko}
\address{Department of Mathematics, University of Notre Dame}
\email{starchenko.1@nd.edu}
\author{Jinhe Ye}
\address{Institut de Math\'ematiques de Jussieu-Paris Rive Gauche}
\email{jinhe.ye@imj-prg.fr}
\date{}
\begin{document}
\begin{abstract}
  We consider \(G\), a linear algebraic group defined over \(\kk\), an 
  algebraically closed field. By considering \(\kk\) as an embedded residue 
  field of an algebraically closed valued field \(K\), we can associate to it 
  a compact \(G\)-space \(S^\mu_G(\kk)\), consisting of \(\mu\)-types on 
  \(G\). We showed that for each \(p_\mu\in S^\mu_G(\kk)\), 
  \(\Stab^\mu(p)=\Stab(p_\mu)\) is a solvable infinite algebraic group when 
  \(p_\mu\) is centered at infinity and residually algebraic. Moreover we 
  give a description of the dimension \(\Stab(p_\mu)\) in terms of dimension 
  of \(p\).
\end{abstract}
\maketitle
\author
\section{Introduction}
Let \(G\) be a group definable in an o-minimal theory, and let 
\(\gamma:(a,b)\ra{}G\) be a definable curve which is unbounded, in the sense 
that the limit at \(b\) does not exist.  In \cite{Kobi}, it was shown that 
one can associate to this datum a definable one-dimensional torsion-free 
group \(H_\gamma\subseteq{}G\), that can be viewed as the ``stabilizer of 
\(\gamma\) at \(\infty\)''. The group \(H_\gamma\) is called the 
\emph{Peterzil--Steinhorn subgroup} associated to \(\gamma\). For example, 
when \(G\) is a Cartesian power of the additive group, \(H_\gamma\) is the 
linear subspace whose translate is the asymptote of \(\gamma\) at \(\infty\).

Assume now that \(G\) is an affine algebraic group over the complex numbers, 
and \(X\) is an algebraic curve embedded in \(G\). If we view \(\mathbb{C}\) 
as the algebraic closure of a real closed field \(\mathcal{K}\), the set of 
complex points of \(X\) can be viewed as the set of \(\mathcal{K}\)-points of 
an \(\mathcal{K}\)-definable set \(X^{an}\) in the o-minimal structure 
\(\mathcal{K}\). This set is unbounded, and we may therefore choose an 
unbounded curve \(\gamma\) inside \(X^{an}\), and consider the corresponding 
PS-group \(H_\gamma\). Taking its Zariski closure, we obtain an algebraic 
subgroup \(G_\gamma\) of \(G\), of (algebraic) dimension \(1\).

It is natural to ask, to which extent does the subgroup \(G_\gamma\) depend 
on the non-algebraic data involved, namely, the dependence on the real closed field \(\mathcal{R}\) of choice and the curve 
\(\gamma\)? And if it does not depend on the above, can the construction be 
described in a purely algebraic manner? We first note that a choice of 
\(\gamma:(a,b)\ra{}X^{an}\) determines an additional \emph{algebraic} datum: 
the curve \(X\) (which we may assume to be smooth) has a canonical 
compactification \(\ti{X}\), its projective model, which is obtained from 
\(X\) by adding finitely many points. Viewing \(\gamma\) as taking values in 
\(\ti{X}^{an}\) rather than \(X^{an}\), the limit of \(\gamma\) at \(b\) will 
be precisely one of these points, and curves \(\gamma\) corresponding to 
different such points definitely might give rise to different subgroups 
\(G_\gamma\). Hence, any hope of providing an algebraic construction of 
\(G_\gamma\) should take into account the choice of such a point at infinity.

\subsection*{Main results}
The main result of this paper (stated below) provides an algebraic 
construction as expected above, once the additional datum of a limit point is 
chosen:

\begin{thm}\label{thm:main1}
  Let \(\kk\) be an algebraically closed field and \(G\) be a linear 
  algebraic group over \(\kk\), and let \(X\subseteq{}G\) be an irreducible 
  curve over \(\kk\). Then there are finitely many \(1\)-dimensional linear 
  subgroups of \(G\), naturally associated to points at infinity of a smooth 
  projective model of the curve \(X\). In fact, they are the 
  \(\mu\)-stabilizers of those points at infinity, as explained later.
\end{thm}

A more precise version is given in Theorem~\ref{thm:main1a}, and the notion 
of the \(\mu\)-stabilizers will be introduced later. The main result of the 
paper, Theorem~\ref{thm:main}, includes a generalization of the above to 
higher dimensions, and some analysis of the structure of the resulting group.  
These results are obtained by viewing \(G\) as a definable group in ACVF, the 
theory of algebraically closed valued fields, and applying some results 
from~\cite{HL}. The relation of the \(\mu\)-stabilizers to the original 
construction of Peterzil--Steinhorn is explained in 
Remark~\ref{rmk:psgroups}.

To state the main result, we need to introduce some additional terminology.  
The subgroups we are interested in were introduced in an abstract setup 
in~\cite{Ser}. There, one considers (suitably defined) definable topological 
groups. To such a group \(G\), one associates an ``infinitesimal subgroup'' 
\(\mu\), the intersection of all definable neighbourhoods of the identity. If 
\(P\) is a (partial) type on \(G\), the set \(\mu{}\cdot P\) can be viewed 
geometrically as a tube around \(P\), and the \(\mu\)-stabilizer 
\(\Stab^\mu(P)\) of \(P\) is defined to be the stabilizer of this set.

In the o-minimal context, the datum of a curve \(\gamma\) as above determines 
a ``type at infinity'' \(p_\gamma\), and it is easy to see that the PS-group 
\(H_\gamma\) depends only on this type. It is shown in~\cite{Ser} that 
\(H_\gamma\) is precisely the \(\mu\)-stabilizer of \(p_\gamma\). Similarly, 
every closed point of the projective model of a smooth curve \(X\) determines 
an ACVF type on \(X\), and the associated group is defined as the 
\(\mu\)-stabilizer of this type. To see that the definition is reasonable, it 
is shown that the resulting group is \(1\)-dimensional. Furthermore, it is 
contained in the (algebraic) stabilizer of the corresponding point in every 
equivariant compactification of \(G\) (Remark~\ref{rmk:comp}).

The definition of \(\mu\)-stabilizer makes sense for types of higher Zariski 
dimensions as well. However, two types of different Zariski dimension might 
have the same tube (Example~\ref{ex:reduced}), so the dimension comparison is 
not straightforward. We say that a type is \(\mu\)-reduced if it is of 
minimal dimension among all types with a given tube. With this terminology, 
we have the following generalization of Theorem~\ref{thm:main1}:

\begin{thm}[Main theorem]\label{thm:main}
  Let \(G\) be a linear algebraic group defined over \(\kk\) and \(p\) be a 
  residually algebraic type. If \(p\) is centered at infinity, then 
  \(\Stab^\mu(p)\) is infinite. Furthermore, if \(p\) is \(\mu\)-reduced, then 
  \(\dim(\Stab^\mu(p))=\dim{}p\). And for each type \(p\), 
  \(\Stab^\mu(p)\) is a solvable linear algebraic group.
\end{thm}

Here the term ``centered at infinity'' should be understood as 
``unbounded'' in the o-minimal counterpart. Do note that one cannot hope for 
the group to be torsion-free, as in the result on PS-groups, since the 
underlying field may have positive characteristic. The three parts of the 
theorem are proved as Cor~\ref{cor:main1}, in \S\ref{sec:dim} and as 
Theorem~\ref{thm:main3}, respectively.
%




\subsection*{Structure of the paper}
The structure of the paper is as follows: In \S\ref{sec:prelim} we review 
some definitions and results related to group actions, and provide an 
alternative approach to \(\mu\)-stabilizers. In \S\ref{sec:onedim} we 
consider the one-dimensional case of Theorem~\ref{thm:main1}. Though formally 
included in the general case, the situation is considerably simpler in this 
case, and sheds light on the more complicated general case. Then in 
\S\ref{sec:general} we deal with the general case.
\subsection*{Acknowledgements}
We would like to thank Antoine Ducros, Fran\c cois Loeser and Kobi Peterzil 
for their helpful conversations. M.~K. was partially supported by the Israel 
Science foundation (Grant no. 1382/15).  J.~Y. was partially supported by NSF 
research grant DMS1500671 and Fran\c cois Loeser's Institut Universitaire de 
France grant.

\section{\texorpdfstring{\(\mu\)}{mu}-stabilizers over ACF}\label{sec:prelim}
Let \(\kk\) be an arbitrary algebraically closed field , and a linear 
algebraic group \(G\) defined over \(\kk\). In this section, we develop the 
theory of \(\mu\)-types and their stabilizers in this context following 
\cite{Ser}. Before going into \(\mu\)-types, let's begin with some generality 
on definable group actions.
\subsection{Definable group actions}\label{ss:gactions}
Let us start by recalling some general facts about stabilizers of definable 
types in an arbitrary complete theory \(T\), following \cite[Section 2]{Ser}. 
We fix a monster model \(\Uu\) of \(T\) and all the models of \(T\) we 
consider will be elementary submodels of \(\Uu\).

Let \(\XX\) be a definable set (over \(0\)), and let \(A\) be a small set of 
parameters. We use \(\LL_\XX(A)\) to denote the set of formulas \(\psi\) over 
\(A\) such that \(\psi(x)\Rightarrow{}x\in{}\XX\). Such formulas will 
occasionally be called \deft{\(\XX\)-formulas}.  And by a \deft{(partial) 
\(\XX\)-type over \(A\)}, we mean a consistent collection of formulas in 
\(\LL_\XX(A)\).

We fix \(\HH\) to be a definable group with a definable action on \(\XX\).  
For an \(\HH\)-formula \(\phi(x)\)  and \(\XX\)-formula \(\psi(y)\), let 
\((\phi\cdot \psi)(z)\) be
\[
\exists x \, \exists y \, \phi(x)\wedge \psi(y) \wedge z=x\cdot y.
\]
And for a partial \(\XX\)-type \(p\), 
\(\phi\cdot{}p=\St{\phi\cdot\psi}{\psi\in{}p}\).

By a \deft{definable \(\XX\)-type over \(A\)}, we mean a \(\XX\)-type over 
\(A\) such that for any formula \(\phi(x,y)\), \(\{a \in A: \phi(x,a)\in 
p\}=\{a \in A: d_p\phi(a)\}\) for some formula \(d_p\phi\) over \(A\). Note 
in the above definition, \(a\) can be tuples in \(A\).

Let \(\Mm\) a model of \(T\) such that \(A=\YY(\Mm)\) for some 
\(\Mm\)-definable set \(\YY\). Then any two such definitions \(d_p\phi\) will 
be equivalent. Moreover, \(p\) can be extended to a unique type  \(p|\Ll\) 
over \(\YY(\Ll)\) determined by \(\St{\phi(x,c)}{c\in{}d_p\phi(\YY(\Ll))}\), 
for any \(\Ll\) such that \(\Mm\preceq \Ll\).

\begin{conv}\label{con}
  For the remainder of the paper, we will assume that we are working over a 
  set of parameters \(A\) such that \(A=\YY(\Mm)\) for some model \(\Mm\), 
  and we assume further that \(\HH\) and its action on \(\XX\) is defined 
  over \(A\).  We assume further that \(\HH \subseteq \YY^n\) for some 
  Cartesian product of \(\YY\).
\end{conv}

\begin{Def}
  Let \(p\) be a definable partial \(\XX\)-type over \(A\) as in 
  Convention~\ref{con}. We define
  \[
    \Stab(p)(\Mm)=\St{h\in\HH(\Mm)}
    {\text{ For any }\phi\in\LL_\XX(A)\,\,p\models h\cdot\phi\Leftrightarrow 
    p\models\phi}
  \]
  where \(h\cdot \phi\) stands for \((x=h)\cdot \phi\). We will occasionally denote \(\Stab(p)(\Mm)\) by \(\Stab(p)(A)\).
\end{Def}

The following is~\cite[Proposition~2.13]{Ser}.
\begin{fact}\label{fact:progroup}
  Let \(\HH\) be a definable group with a definable action on \(\XX\) and 
  assume we are in the setting of Convention~\ref{con}. Let \(p\) be a 
  partial definable \(\XX\)-type over \(A\). Then \(\Stab(p)\) is a 
  \(A\)-type-definable subgroup of \(\HH\) in the following sense: there is a 
  small system \(\HH_\alpha\) of \(A\)-definable subgroups of \(\HH\) such 
  that for every elementary extension \(\Mm\preceq \Ll\), for 
  \(a\in\HH(\Ll)\), we have that \(a \in \Stab(p|\Ll)(\Ll)\)
  if and only if \(a\in\HH_\alpha(\Ll)\) for all \(\alpha\).
\end{fact}
With the language set up, we will now look at the setting to talk about 
\(\mu\)-types as in \cite{Ser} over algebraically closed fields.

\subsection{\texorpdfstring{\(\mu\)}{mu}-stabilizers over ACF}
Let \(\kk\) be an algebraically closed field. The theory of algebraically 
closed fields is not rich enough to have a good notion of infinitesimal 
subgroups as in \cite{Ser}.  Hence, it is natural to work with the theory 
\(T_{loc}\) as introduced in \cite[Section~6]{HK3}.  The language for 
\(T_{loc}\) has two sorts, a sort \(\VF\) for the valued field, and a sort 
\(\Gamma\) for the value group. The sort \(\VF\) is has a unary predicate 
\(\RES\) for an embedded copy of the residue field, which can be viewed as an 
additional sort. Thus we have function symbols \(\res:\VF^2\to\RES\) and 
\(\val:\VF\to\Gamma\). The theory \(T_{loc}\) asserts that the \(\VF\) sort 
is an algebraically closed valued field, \(\RES\) is a subfield, \(\val\) is 
a valuation map, \(\res(x,y)=\res(x/y)\), the residue of \(x/y\) if 
\(\val(x)\geq\val(y)\) and \(0\) otherwise,
with \(\res(c,1)=c\) for \(c\in\RES\). For notational simplicity, we will use 
\(\res(x)\) to denote \(\res(x,1)\) for \(x\in\OO\), the valuation ring.

We further assume that we have constants for the elements of the field 
\(\kk\) in \(\RES\). Thus, models of \(T_{loc}\) are algebraically closed 
valued fields with embedded residue field extending \(\kk\).

\begin{fact}[{\cite[Lemma 6.3]{HK3}}]\label{fact:QE}
  \(T_{loc}\) admits quantifier elimination in this language. The sorts 
  \(\Gamma\) and \(\RES\) are stably embedded and orthogonal to each other.
  The induced structures on \(\Gamma\) and \(\RES\) are of divisible ordered 
  abelian groups algebraically closed fields, respectively.
\end{fact}

\begin{rmk}
  In the paper \cite{HK3}, a constant symbol \(1\) in the \(\Gamma\)-sort for 
  some positive element was included. But the proof of quantifier elimination 
  does not rely on the constant.
\end{rmk}

In some cases, we will work in the reduct of \(T_{loc}\) in the 3-sorted 
language \(\LL_\val\). \(\LL_\val\) consists of the valued field sort 
\(\VF\), the value group sort \(\Gamma\) and the residue field sort \(\RES\) 
and maps \(\val:\VF\to\Gamma\), \(\res:\VF\to\RES\). In this language, we 
have constants for \(\kk\) in both \(\VF\)-sort and \(\RES\)-sort. The 
induced theory on the reduct is \(\ACVF_\kk\), the theory of algebraically 
closed valued fields with constants for \(\kk\), which admits quantifier 
elimination in \(\LL_\val\). Since this is a reduct, we will freely view 
formulas in \(\LL_\val\) as definable sets in \(T_{loc}\). The main point of 
working with this reduct is the application of topological results 
from~\cite{HL} available for this restricted class of definable sets.

Recall that we are given a linear algebraic group \(G\) defined over \(\kk\).
In our context, this group determines a number of distinct definable groups: 
The definable group \(\GV\) of \(\VF\)-points of \(G\), and the definable 
subgroups \(\GO\) and \(\GR\) of \(\OO\)- and \(\RES\)-points, respectively.  
The (pointwise) residue map \(\res:\GO\rightarrow\GR\) is a definable group 
theoretic retraction (since \(\GV\) is over \(\kk\)), whose kernel we denote 
by \(\mu\). Note \(\mu\) is definable over \(\kk\) as well.  Geometrically, 
\(\mu\) can be viewed as an (infinitesimal) neighbourhood of the identity in 
\(\GV\).

Let \(P\) be an algebraic variety over \(\kk\) on which \(G\) acts (Our main 
example will be \(P=G\), with \(G\) acting on itself by left multiplication, 
but occasionally we will need the more general setup.) Then \(P\) determines 
a definable set \(\PP\) in \(\VF\), and a definable action of \(\GV\) on 
\(\PP\), which restricts to an action of \(\GR\leq\GV\). We are therefore in 
the setting of Convention~\ref{con}, where we set \(T=T_{loc}\),
\(\XX=\PP\), \(\YY=\RES\) and \(\HH=\GR\). In fact, one can take \(\kk=\RES(\Mm)\), where 
\(\Mm\) is a  field of Hahn series with coefficients in \(\kk\).a

\begin{Def}
  Let \(P\) and \(G\), and the associated terminology be as above. We denote by \(S_\PP(\kk)\) the space of complete 
  \(\LL_\val\)-\(\PP\)-types over \(\kk\). For \(p\in{}S_\PP(\kk)\), the 
  \deft{\(\mu\)-stabilizer} \(\Stab^\mu(p)\) of \(p\) is 
  \(\Stab_{\GR}(\mu\cdot{}p)\).
\end{Def}

Note that a complete \(\LL_\val\)-\(\PP\)-type might be a partial type in 
\(T_{loc}\). By quantifier elimination, such types correspond to pairs of the 
form \((\ZZ,v)\), where \(\ZZ\) is an irreducible closed subvariety of 
\(\PP\) and \(v\) is a valuation on the function field of \(\ZZ\) which is 
trivial on \(\kk\). By Fact~\ref{fact:QE}, the \(\RES\)-sort is stably 
embedded as an algebraically closed field. In particular, each 
\(p\in{}S_\PP(\kk)\) is definable over \(\kk\) in \(\LL_\val\).

\begin{prop}
  Let \(p\in S_\PP(\kk)\). Then \(\mu\cdot{}p\) is a definable partial type 
  over \(\kk\).
\end{prop}
\begin{proof}
  Let \(\ZZ\) be a \(\LL_\val\)-definable set over \(\kk\). Then 
  \(\mu\cdot{}p\models\ZZ\) if and only if  
  \[
  p(x)\models\forall\epsilon\in\mu(\epsilon\cdot x\in \ZZ)
  \]
  The latter condition is \(\LL_\val\)-definable over \(\kk\), hence the 
  result follows from the definability of \(p\).
\end{proof}

By Fact~\ref{fact:progroup} and the above discussion, \(\Stab^\mu(p)\) is 
given by an intersection of \(\LL_\val\)-definable subgroups of \(\GR\).  
However, \(\GR\) has the descending chain condition on subgroups. Hence by Fact 
\ref{fact:QE} we have:
\begin{cor}
  Let \(p\) be an \(\LL_{val}\)-complete \(G\)-type over \(\kk\).  Then the 
  \(\mu\)-stabilizer \(\Stab^\mu(p)\) of \(p\) is a \(\kk\)-definable 
  subgroup of \(\GR\), in the sense that there is a \(\kk\)-definable 
  subgroup \(\mathbf{H}\) of \(\GR\), such that 
  \(\Stab^\mu(p|\Ll)(\Ll)=\mathbf{H}(\Ll)\) for any model \(\Ll \succeq 
  \Mm\).
\end{cor}

For \(p\) and \(q\in{}S_\PP(\kk)\), define \(p\sim{}q\) if 
\(\mu\cdot{}p=\mu\cdot{}q\). It is easy to check that 
\(\mu\cdot{}p=\mu\cdot{}q\) iff in a monster model \(\Uu\), there are 
\(a\models{}p\), \(b\models{}q\) and \(\epsilon\in\mu\) such that 
\(\epsilon\cdot{}a=b\).

We denote by \(S^\mu_\PP(\kk)\) the quotient by this equivalence relation and 
for each \(p\in{}S_\PP(\kk)\), we denote by \(p_\mu\) its equivalence class.  
Since \(\mu\) is normal in \(\GV(\OO)\), the \(\GR(\kk)\) action on 
\(S_\PP(\kk)\) given in Subsection 2.1 respects the equivalence relation.  
Hence \(\GR(\kk)\) acts on \(S^\mu_\PP(\kk)\) , and 
\(\Stab^\mu(p)=\Stab(p_\mu)\), where the right hand side is by considering 
\(\GR(\kk)\) acting on \(S^\mu_\PP(\kk)\).

Lastly, we note that stabilizers for types in the same orbit are conjugate:

\begin{lem}\label{lem:conj}
  Let \(g\in\GR(\kk)\) be such that \(g\cdotp{}p_\mu=q_\mu\). Then 
  \(\Stab^\mu(q)=g\Stab^\mu(p)g^{-1}\).
\end{lem}

We finish this subsection with two basic examples. These can be compared to 
the computations in the o-minimal setting.

\begin{eg}\label{ex:basic1}
  Let \(G=\SL_{2,\kk}\). Let
  \[
  X_1=\St{\begin{pmatrix}
    x & 1 \\
    0 & x^{-1}
  \end{pmatrix}}{x\in\Gm},
  \]
  a closed subvariety of \(G\). Since this subvariety is isomorphic (as an 
  algebraic variety) to \(\Gm\), the definable subset given by \(\val(x)<0\)  
  isolates a complete type \(p\) on \(\GV\) (in the language \(\LL_\val\) 
  over \(\kk\)).

  We claim that the left $\mu$-stabilizer \(H\) of $p$ is the subgroup
  \[
    G_1 =\left\{
    \begin{pmatrix}
      1 & a \\
      0 & 1
    \end{pmatrix}
    \right\}.
  \]

  Indeed, let \(g\in\RES\) and choose $\alpha\in\VF$ with $\val(\alpha)<0$.  
  Then
  \[
  \begin{pmatrix}
    1 & g \\ 0  & 1 
  \end{pmatrix}
  \begin{pmatrix}
    \alpha & 1 \\ 0 & \alpha^{-1} 
  \end{pmatrix}=
  \begin{pmatrix}
    \alpha & 1+g\alpha^{-1} \\
    0 & \alpha^{-1}
  \end{pmatrix}=
  \begin{pmatrix}
    1+\varepsilon & 0 \\
    0 & (1+\varepsilon)^{-1}
  \end{pmatrix}
  \begin{pmatrix}
    \beta & 1 \\ 0 & \beta^{-1}
  \end{pmatrix},\]
  where $\varepsilon =g \alpha^{-1}$ and
  $\beta=(1+\varepsilon)^{-1}\alpha$.
  Since \(\val(\varepsilon)>0\), we have
  \(\Mat{1+\varepsilon}{0}{0}{(1+\varepsilon)^{-1}}\in\mu\) and 
  \(\Mat{\beta}{1}{0}{\beta^{-1}}\models{}p\).

  Thus \(G_1\subseteq H\). By (the easy part of) Theorem~\ref{thm:main1a}, 
  \(H\) is \(1\)-dimensional, so \(G_1\) is the connected component of \(H\), 
  and \(H\) has the form
  \[
    H=\left\{
    \begin{pmatrix}
      \xi & a \\
      0 & \xi^{-1}
    \end{pmatrix} : \xi^n=1
    \right\}.
  \]
  for some \(n\). By a similar computation, such an element will take 
  \(\mutube{p}\) to \(\mutube{p_\xi}\) where \(p_\xi\) is the type of 
  elements
  \(\Mat{x}{\xi}{0}{y}\) with \(\val(x)<0\), so we must have \(\xi=1\) and 
  \(G_1=H\).

  This example will be considered again in~\ref{ex:basic1a}
\end{eg}

\begin{eg}\label{ex:basic2}
  Similarly to Example~\ref{ex:basic1}, we now consider the closed subvariety
  \[
  X_2=\St{\begin{pmatrix}
    x & 0 \\
    1 & x^{-1}
  \end{pmatrix}}{x\in\Gm}\subseteq G=\SL_2,
  \]
  and let \(q\) the type on \(\GV\) determined inside it by the condition 
  \(\val(x)<0\).

  We claim that the left $\mu$-stabilizer \(H\) of $q$ is now the subgroup
  \(G_2=\St{\Mat{a}{0}{0}{a^{-1}}}{a\in\Gm}\). The computation is similar: 
  For $g\in\Gm(\RES)$ and $\alpha\in\VF$ with $\val(\alpha)<0$ we have
\[
  \begin{pmatrix}
    g & 0 \\ 0  & g^{-1} 
  \end{pmatrix}
  \begin{pmatrix}
    \alpha & 0 \\ 1 & \alpha^{-1} 
  \end{pmatrix}=
  \begin{pmatrix}
    g\alpha & 0 \\
    g^{-1} & (g\alpha)^{-1}
  \end{pmatrix}=
  \begin{pmatrix}
    1 & 0 \\
    \varepsilon & 1
  \end{pmatrix}
  \begin{pmatrix}
    \beta & 0 \\ 1 & \beta^{-1}
  \end{pmatrix},\]
  where $\beta=g\alpha$ and $\varepsilon =\frac{g^{-1}-1}{g\alpha}$.
  
  Since $\val(\varepsilon)>0$ and $\val(\beta)<0$, we obtain 
  \(G_2\subseteq{}H\).  The only one-dimensional algebraic subgroup of \(G\) 
  that properly contains \(G_2\) is its normalizer, which contains the 
  element \(w=\Mat{0}{1}{-1}{0}\). This elements sends \(q\) to the type of 
  elements \(\Mat{1}{x^{-1}}{-x}{0}\) with \(\val(x)<0\), so cannot be in 
  \(H\). Hence the \(\mu\)-stabilizer is \(G_2\) in this case.
\end{eg}

\subsection{A different view on \texorpdfstring{\(\mu\)}{mu}-stabilizers}
Instead of viewing the \(\mu\)-stabilizers syntactically as in the previous 
section, we have some concrete constructions to realize them in the monster 
model as well. In this section, we describe the construction, following the 
same idea of~\cite[\S2.4]{Ser}. We work in a fixed monster model of 
\(T_{loc}\), \(\Uu\), and identify definable sets and (partial) types with 
their realisations in \(\Uu\). From now on, we restrict our attention to the 
case \(\PP=\GV\), unless mentioned otherwise.

\begin{Def}
  For \(p\in{}S_G(\kk)\) we use \(\GR_p\) to denote the following set 
  \(((\mutube{p})\cdotp(\mutube{p})^{-1})\cap\GR\).
\end{Def}

\begin{prop}\label{prp:Gp}
  Let \(a\in\mutube{p}\). The following are equivalent for an element 
  \(b\in\GR(\kk)\):
  \begin{enumerate}
    \item\label{prp:Gp:Rf} \(b\in\mutube{p}\cdot{}a^{-1}\cap\GR\)
    \item\label{prp:Gp:Of} \(b=\res(a_1a^{-1})\) for some \(a_1\models{}p\) 
      for which \(a_1a^{-1}\in\GO\)
    \item\label{prp:Gp:Rg} \(b\in\GR_p\)
    \item\label{prp:Gp:Og} \(b=\res(a_1{a_2}^{-1})\) for some 
      \(a_1,a_2\models{}p\) for which \(a_1{a_2}^{-1}\in\GO\)
  \end{enumerate}
  Hence, \(\GR_p(\kk)=\mutube{p}\cdotp{}a^{-1}\cap\GR(\kk)\)
\end{prop}
\begin{proof}
  The equivalence of~\eqref{prp:Gp:Rf} and~\eqref{prp:Gp:Of} follows directly 
  from the definitions, and likewise for~\eqref{prp:Gp:Rg} 
  and~\eqref{prp:Gp:Og}. Hence we need to show~\eqref{prp:Gp:Og} 
  implies~\eqref{prp:Gp:Of}. Assume \(b=\res(a_1a_2^{-1})\). Since \(a\) and 
  \(a_2\) satisfy the same type over \(\kk\), and \(\RES\) is stably 
  embedded and stable, there is an automorphism \(\tau\) over \(\kk\) such that 
  \(\tau(a_2)=a\). Then \(b=\tau(b)=\res(\tau(a_1)a^{-1})\), with 
  \(\tau(a_1)\) also satisfying \(p\), showing~\eqref{prp:Gp:Of} for \(b\).
\end{proof}

We now have the following description of \(\GR_p\).

\begin{cor}\label{def}
  \(\GR_p(\kk)=\Stab^\mu(p)(\kk)\)
\end{cor}
\begin{proof}
  Assume \(g\in\GR\) stabilizes \(\mutube{p}\). Then for any 
  \(a\in\mutube{p}\), \(g\cdot a\in\mutube{p}\), hence \(g\in\mutube{p}a^{-1}\), so 
  is in \(\GR_p\). Conversely, if \(g\in\GR_p\) and \(a\in\mutube{p}\), 
  writing \(g=a_1a^{-1}\) as above we obtain \(g\cdot a\in\mutube{p}\).
\end{proof}

\begin{rmk}\label{rmk:def}
  We would like to have Cor.~\ref{def} to hold for \(\RES\)-points, instead 
  of just \(\kk\)-points.  This is not automatic, since \(\GR_p\) is not, 
  a-priori, a definable set. However, in the special case as in Theorem 
  \ref{thm:main1}, it is indeed the case. 

\end{rmk}

\subsection{\texorpdfstring{\(\mu\)}{mu}-reduced types}
For \(p\in{}S_G(\kk)\) we denote by \(\dim(p)\) the dimension of its Zariski 
closure in \(\GV\) over \(\kk\). Most of the rest of this paper is devoted to 
comparing this dimension to the dimension of \(\Stab^\mu(p)\). We first note 
that if \(\XX\) is a variety over a valued field \(L\) whose valuation ring 
is  \(\OO\), then the Zariski dimension of \(\res(\XX\cap\OO)\) is at most 
the dimension of \(\XX\) (this follows for example from~\StackP[Lemma]{00QK}, 
by choosing a model of \(\XX\) over \(\OO\)).  Applying this observation to 
\(\XX=\YY{}a^{-1}\), where \(\YY\) is a variety containing \(p\), we obtain:

\begin{prop}\label{Dim}
  For any \(p\in{}S_G(\kk)\), 
  \(\dim(\Stab^\mu(p))=\dim\GR_p(\kk)\leq\dim(p)\), where \(\dim\) means the 
  Krull dimension in \(\Stab^\mu(p)\) and \(\GR_p(\kk)\) respectively.  
  \(\dim(p)\) is the minimal \(\VF\)-dimension of the formulas \(\varphi \in 
  p\)
\end{prop}

In general, the above bound will not be sharp, since types of different 
dimensions may have the same \(\mu\)-type:

\begin{eg}\label{ex:reduced}
  Let \(G=\AAA^2\) as an additive group. Let \(K\) be a large enough Hahn 
  series in variable \(t\) over \(\kk\). Let 
  \(p=tp((t^{-1},t^{-1}+t^r)/\kk)\) where \(r>0,r\notin\mathbb{Q}\), and 
  \(tp\) denotes the \(\LL_\val\)-type. Then \(\dim(p)=2\), since 
  \(t^{-1}+t^r\) is transcendental over \(t^{-1}\). But 
  \(\mutube{p}=\mutube{q}\), where \(q=tp((t^{-1},t^{-1})/\kk)\), since 
  \((t^{-1},t^{-1}+t^r)\) and \((t^{-1},t^{-1})\) differ by
  \((0,-t^r)\in\mu\), so \(\dim(\GR_p)\le{}1\) (as we will later see, they 
  are in fact equal).  Furthermore, when \(\mathrm{Char}(k)=p>0\), we see 
  that \(\Stab^\mu(p)\) is not torsion-free.
\end{eg}

This observation motivates the following definition.

\begin{Def}
  For \(p\in{}S_G(\kk)\), we say that \(p\) is \deft{\(\mu\)-reduced} if 
  \(p\) is a type of minimal dimension in \(p_\mu\). An element \(a\in\GV\) 
  is \(\mu\)-reduced over \(\kk\) if \(a\models p\) for some \(\mu\)-reduced \(p\).
\end{Def}
\subsection{Bounded types}
In this sub-section, we revert to working with a general \(G\)-variety 
\(\PP\). We recall the following definition (e.g., from~\cite[\S4.2]{HL}):
\begin{Def}
  Let \(\VV\) be an affine variety, viewed as a definable set in ACVF, and 
  let \(\XX\subseteq\VV\) be a \(\LL_\val\)-definable subset. We say that \(\XX\) is 
  \deft{bounded} if for every regular function \(f\) on \(\VV\) there is 
  \(\gamma\in\VG\) such that \(\val(f(\XX))\geq\gamma\).
  
  For a general variety \(\VV\), a subset \(\XX\subseteq\VV\) is bounded if 
  it is covered by bounded subsets of an affine cover.
  
  A partial type \(p\) in \(\VV\) is \deft{bounded} if \(p\models\XX\) for 
  some bounded \(\XX\subseteq\VV\). A type in \(\VV\) is said to be 
  \deft{centered at infinity} if it is not bounded.
\end{Def}

Note that the property of a definable set to be bounded depends on the 
ambient variety (for example, \(\AAA^1\) is bounded as a subset of 
\(\mathbb{P}^1\), but not as a subset of \(\AAA^1\)). However, if 
\(\VV\) is a closed subvariety of \(\WW\), then \(\XX\subseteq\VV\) is 
bounded in \(\VV\) if and only if it is bounded in \(\WW\). Also, it suffices 
to check the conditions for generators of the regular functions. In 
particular, a subset \(\XX\) of a closed subvariety of \(\AAA^n\) is bounded 
if and only if \(\val(\XX)\ge\gamma\) for some \(\gamma\).

Over \(\kk\), we have in our situation the following:

\begin{prop}
  A \(\kk\)-definable set \(\XX\subseteq\VV\) is bounded if and only if it is 
  contained in \(\VV(\OO)\)
\end{prop}
\begin{proof}
  By definition, it suffices to prove the statement for \(\VV\) affine, and 
  by the remarks above, for \(\VV=\AAA^n\).

  If \(\XX\subseteq\OO^n\) we may take \(\gamma=0\) in the definition.  
  Conversely. We may assume \(n=1\) by projecting. If \(a\in{}\XX\setminus\OO\) 
  then \(\gamma=\val(a)<0\) has the same type as any other negative value 
  \(\gamma'\), so there is an automorphism of \(\VG\) taking \(\gamma\) to 
  \(\gamma'\), and since \(\VG\) is stably embedded and \(\VG\) and \(\RES\) 
  are orthogonal, it extends to an automorphism over \(\kk\) that takes \(a\) 
  to \(a'\in{}\XX\), with \(\val(a')=\gamma'\).  Thus, \(\XX\) is unbounded.
\end{proof}

Let \(p\) be a bounded type on \(\PP\), a variety endowed with an action of 
\(\GV\). A realization \(a\) of \(p\) is then an \(\OO\)-point of \(\PP\), 
and so determines a point \(\ol{a}\) of \(\PP\) in the residue field. The 
type of \(\ol{a}\) depends only on \(p\) (since it is encoded there), and we 
denote it by \(\ol{p}\). The group \(\GR\) acts on the set of all types in 
\(\ol{\PP}\), the variety \(\PP\) viewed as a definable set in \(\RES\). In 
particular, we may consider the stabilizer of \(\ol{p}\).

\begin{prop}\label{prp:stabres}
  For any bounded type \(p\) on \(\PP\) we have 
  \(\Stab^\mu(p)\leq{}\Stab(\ol{p})\).
\end{prop}
\begin{proof}
  Let \(\ol{a}\) be a realization of \(\ol{p}\), and let \(a\) be a 
  realization of \(p\) whose residue is \(\ol{a}\). Assume that for some 
  \(g\in\GR\) we have \(g\cdot a=\epsilon{}\cdot b\) for some \(\epsilon\in\mu\) and 
  \(b\) realizing \(p\) (so that \(g\in\Stab^\mu(p)\)).  Since all elements 
  involved are in \(\OO\), we may apply the residue map, and obtain 
  \(g\cdot\ol{a}=\ol{b}\).  Since \(b\) realizes \(p\), \(\ol{b}\) realizes 
  \(\ol{p}\). Thus \(g\cdot \ol{p}=\ol{p}\), i.e., \(g\in\Stab(\ol{p})\).
\end{proof}

Returning to the case \(\PP=\GV\), we obtain:

\begin{cor}
  If \(p\) is a bounded type on \(\GV\) such that \(\ol{p}\) is realized in 
  \(\kk\) then its \(\mu\)-stabilizer is trivial.
\end{cor}
\begin{proof}
  In this case \(\ol{p}\) corresponds to a (closed) point of \(G\), hence the 
  stabilizer is trivial.
\end{proof}

Because of the last corollary, we shall concentrate on types centered at 
infinity.

\section{Analyzing the one-dimensional case}\label{sec:onedim}
In this section, we prove the main theorem in dimension~1 
(Theorem~\ref{thm:main1}). We think this section worth including even 
though it follows from the general case, since it is relatively simple, and 
it sheds light on the important idea in proving the general case. The result 
in this section is first proved by Moshe Kamensky and Sergei Starchenko in 
their unpublished notes via the language of places.

\subsection{Points on curves}
Each smooth curve \(\XX\) over \(\kk\) embeds in a unique smooth projective 
one over \(\kk\), its projective model \(\ti{\XX}\). Every closed point \(c\) 
on \(\ti{\XX}\) corresponds to a valuation \(\val_c\) on the function field 
\(\kk(\XX)\), given by the order of vanishing at \(c\). In particular, 
\(\val_c\) is trivial on \(\kk\). The projective model contains a finite 
number of closed points outside of \(\XX\), which we call the points at 
infinity.

In our case, \(\XX\) is an affine curve, embedded as a closed subvariety in a 
fixed affine space \(\AAA^n\). To any point \(c\in\ti{\XX}\) we associate the 
complete type on \(\XX\) determined by
\[
  p_c(a)=\St{\val(f(a))>0}{\val_c(\ol{f})>0}
\]
where \(f\) runs over all elements of the local ring corresponding to 
\(\XX\), and \(\ol{f}\) is the corresponding element in \(\kk(\XX)\).

We would like to describe the types that occur in this way intrinsically, in 
a way that will be helpful later. The condition that \(c\) is a closed point 
corresponds to the following.

\begin{Def}
  An extension of (possibly trivially) valued fields is \deft{residually 
    algebraic} if the corresponding residue field extension is algebraic.  
  For \(L\) a (possibly trivially) valued field, an \(\LL_\val\)-type \(p\) 
  over \(L\) is \deft{residually algebraic} if  for a/every realization 
  \(a\) of \(p\), \(L(a)\) is residually algebraic over \(L\).
\end{Def}

\begin{prop}
  Let \(\XX\) be a smooth curve embedded in \(\AAA^n\) (viewed as a definable 
  set in \(\VF^n\)). A \(\LL_\val\)-type \(p\) over \(\kk\)  on \(\XX\) is 
	residually algebraic if and only if it is of the form \(p_c\) for a closed 
	point \(c\) of \(\ti{\XX}\), the smooth projective model of \(\XX\).  
	Furthermore, \(c\in\ti{\XX}\setminus\XX\) if and only if \(p_c\) is 
	unbounded.
\end{prop}
\begin{proof}
  Let \(p\) be a residually algebraic type on \(\XX\), and let \(a\) be a 
  realization that witnesses this. If \(a\in\kk\), \(p\) corresponds to the 
  \(\kk\)-point \(a\) of \(\XX\), and we are done. Otherwise, \(\kk(a)\) is 
  isomorphic to \(\kk(\XX)\) as a field, and since \(p\) is residually 
  algebraic, the valuation on \(\kk(a)\) is non-trivial. Thus, we obtain a 
  \(\kk\)-point \(c\) of \(\ti{\XX}\) by the discussion earlier in this 
  subsection, and it is clear that the two procedures are inverse to each 
  other. The last statement also follows from the above.
\end{proof}

We have been working with smooth curves, but since we are interested in 
points at infinity, hence the assumption is immaterial, since the singularity of varieties are at least codimension 1, hence varieties are smooth 
at generic points.

\begin{cor}\label{cor:fin}
  For \(\XX\) a curve, there are finitely many residually algebraic types 
  centered at infinity. Moreover, they are isolated by \(\LL_\val\)-formulas 
  over \(\kk\).
\end{cor}
\begin{proof}
  It remains to prove the moreover part of the statement. Since we know that there 
  are only finitely many types on \(\XX\) centered at infinity, call them 
  \(p_1,...,p_m\). Without loss of generality, for each \(i,j\), there will 
  be regular functions \(f_{ij},g_{ij}\) such that 
  \(p_i\models\val(f_{ij})<\val(g_{ij})\) but 
  \(p_j\models\val(f_{ij})\geq\val(g_{ij})\). Hence some Boolean combinations 
  of the formulas above together with the formula \(x\notin\XX(\OO)\) will 
  isolate the above types.
\end{proof}

\begin{rmk}\label{rmk:comp}
  Let \(\XX\) be an affine curve embedded in \(G\), and assume that we are 
  given a \(G\)-equivariant embedding of \(G\) in a \(G\)-variety \(P\).  
  Assume that the closure \(\XX'\) of \(\XX\) in \(P\) includes the point 
  \(c\in\ti{\XX}\). The type \(p_c\) is then bounded in \(P\), and by 
  Prop.~\ref{prp:stabres} the \(\mu\)-stabilizer of \(p_c\) is contained in 
  the stabilizer of the residue type of \(p_c\), which is simply \(c\).  
  Hence, the \(\mu\)-stabilizer of \(c\) is contained in the stabilizer of 
  this point in every equivariant ``compactification'' where the point is 
  realized.

  This fact, along with the dimension equalities for the \(\mu\)-stabilizers 
  justifies viewing the \(\mu\)-stabilizer as a ``canonical stabilizer'' for 
  the corresponding point.
\end{rmk}

\begin{eg}\label{ex:basic1a}
  Example~\ref{ex:basic1} provides an instance of the situation here, with 
  \(\XX=X_1\). Since \(X_1=\Gm\) as an abstract variety, its projective model 
  \(\ti{X_1}\) is \(\PP^1\), so has two additional points, \(0\) and 
  \(\infty\) (with \(0\) the one included in the chart where \(x\) is 
  defined). Hence, the type \(p\) considered there is \(p_\infty\) in our 
  notation, and \(p_0\) corresponds to the type with the roles of \(x\) and 
  \(x^{-1}\) reversed.
  
  We may alter it a bit by considering the image of \(X_1\) in \(\PSL_2\) 
  (the computation remains essentially the same, but we now also have the 
  elements \(\Mat{-1}{a}{0}{-1}\) in the stabilizer). The space \(\PSL_2\) 
  can be compactified equivariantly by mapping it into \(\PPP^3\), viewed as 
  the projective space associated to the space of all linear endomorphisms of 
  \(\AAA^2\). It is easy to compute that under this embedding, the point 
  \(\infty\) maps to the element \(c=[\Mat{1}{0}{0}{0}]\) of \(\PPP^3\). The 
  stabilizer of this point under the action of \(\SL_2\) is the subgroup of 
  upper triangular matrices, so properly contains our \(\mu\)-stabilizer.
\end{eg}

Let \(p\in S_G(\kk)\) be a residually algebraic type of Zariski dimension 
\(1\) inside \(\GV\). There is then a curve \(\XX\) in \(\GV\) containing 
\(p\). We had explained in Proposition~\ref{prp:Gp} and Remark~\ref{rmk:def} 
that \(\Stab^\mu(p)(\kk)=\mutube{p}\cdot{}a^{-1}\cap\GR(\kk)\) for any 
realization \(a\) of \(p\) (this will be shown again for residually algebraic 
types in Cor.~\ref{model}). However, since \(p\) is isolated by 
Corollary~\ref{cor:fin}, we see that 
\(\Stab^\mu(p|\Ll)(\Ll)=\mutube{p}\cdot{}a^{-1}(\Ll)\cap\GR(\Ll)\) for any 
\(\Ll\) extending \(\Mm\) and \(a\). In particular, one can work with a model 
\(\Ll\) of \(T_{loc}\) with \(\RES(\Ll)=\kk\). Working in this model, let 
\(p_i\), \(p_j\) be two types as above. If \(g\in G(\kk)\) satisfies 
\(g\cdot{}p_i\in{}{p_j}_\mu\) (i.e., if \(\mu\cdot g\cdot p_i=\mu\cdot 
p_j\)), then \(\mu\cdot{}p_j\cdot{}a^{-1}=g\cdot\Stab^\mu(p_i)\), for any 
\(a\models{}p_i\).

To complete the proof, we would like to show that this set 
is infinite for some realization \(a\) of \(p\). This amounts to showing that 
\(\mutube{p}\cdot{}a^{-1}\) cannot be covered by a finite number of open 
balls. To do that, we will use topological methods from~\cite{HL}, which we 
review below.

\subsection{Tame topology on definable sets}\label{sec:tame}
We make a slight digression into the tame topology of definable sets in 
\(\ACVF\), as developed in~\cite{HL}. This is an important ingredient in the 
proof of the main result.

The results in this section can be found in \cite{HL}. In this section, the underlying 
theory is \(\ACVF\), and the main motivation is to study the topological 
structure of \(\LL_\val\)-definable sets in the \(\VF\)-sort.

\begin{Def}
  Let \(\VV\) be an algebraic variety over a valued field \(F\), a subset 
  \(\XX\subseteq{}\VV\) is \deft{\(v\)-open} if it is open for the valuative 
  topology.
  
  A subset \(\XX\subseteq\VV\) is \deft{\(g\)-open} if it is a positive 
  Boolean combination of Zariski open, closed sets and sets of the form
  \begin{align*}
     \St{x}{v\circ f(x)>v\circ g(x)}
  \end{align*}
  where \(f\) and \(g\) are regular functions defined on \(\mathbf{U}\), a 
  Zariski open subset of \(\VV\).
  
  If \(\ZZ\subseteq\VV\) is a definable subset of \(\VV\), a subset 
  \(\mathbf{W}\) of \(\ZZ\) is said to be \(v\)(respectively \(g\))-open if 
  \(\mathbf{W}\) is of the form \(\ZZ\cap\YY\), where \(\YY\) is \(v\) 
  (respectively \(g\))-open in \(\VV\).
  
  The complement of a \(v\)(respectively \(g\))-open is called 
  \(v\)(respectively \(g\))-closed.  We say \(\XX\) is \(v{+}g\)-open 
  (respectively \(v{+}g\)-closed) if it is both \(v\)-open and \(g\)-open 
  (respectively both \(v\)-closed and \(g\)-closed).
\end{Def}

Note that the \(v{+}g\)-opens does not form a topology, as it is not even 
closed under arbitrary union. However, it is still makes sense to talk about 
connectedness in this setting:

\begin{Def}
  Let \(\XX\) be a definable subset of \(\VV\), an algebraic variety. We say 
  that \(\XX\) is \deft{definably connected} if \(\XX\) cannot be written as 
  a disjoint union of two non-empty \(v{+}g\)-open subsets of \(\XX\).
  
  We say that \(\XX\) has \deft{finitely many definably connected components} 
  if \(\XX\) can be written as a finite disjoint union of \(v{+}g\)-clopen 
  definably connected subsets.
\end{Def}

\begin{Def}
  Let \(f:\VV\rightarrow\WW\) be a definable function from \(\VV\) to 
  \(\WW\), we say \(f\) is \deft{\(v\)-continuous} if \(f^{-1}(\XX)\) is 
  \(v\)-open for \(\XX\) a \(v\)-open subset of \(\WW\), and we define 
  \deft{\(g\)-continuous} functions similarly. We say \(f\) is 
  \(v{+}g\)-continuous if \(f\) is both \(v\)-continuous and 
  \(g\)-continuous.
\end{Def}

\begin{prop}[Hrushovski, Loeser]\label{prp:defcon}
  If \(f\) is \(v{+}g\)-continuous and \(\XX\) is a definably connected and 
  \(f\) is defined on \(\XX\), then \(f(\XX)\) is definably connected.

  If \(\VV\) is an geometrically/absolutely irreducible variety, then \(\VV\) 
  is definably connected.
\end{prop}

The following is an easy corollary of \cite[Theorem 11.1.1]{HL}. 

\begin{thm}[Hrushovski, Loeser]\label{comp}
  Given a definable subset \(\XX\subseteq \VV\), where \(\VV\) is some quasi-projective variety, \(\XX\) has finitely many definably connected components.
\end{thm}
We also have the following.
\begin{thm}\label{thm:boundedvar}
  Let \(\VV\subseteq\AAA^n\) be a closed subvariety. \(\VV\) is bounded 
  iff \(\VV\) is zero dimensional.
\end{thm}
\begin{proof}
  If \(\VV\) is bounded, then \(\VV\) will be definably compact as in 
  Section~4 of \cite{HL}. This implies \(\VV\) is proper by 
  \cite[Proposition~4.2.30]{HL}, hence \(\VV\) is zero dimensional. The 
  converse is clear.
\end{proof}
We may now prove the following more precise version of 
Theorem~\ref{thm:main1} (the case of curves).

\begin{thm}\label{thm:main1a}
  Let \(p\in S_G(\kk)\) be a residually algebraic type, centered at infinity 
  with \(\dim(p)=1\). Then \(\dim(\Stab^\mu(p))=1\).
\end{thm}
\begin{proof}
  By Prop.~\ref{Dim} and the discussion preceding Section~\ref{sec:tame}, it 
  suffices to show that the \(\mu\)-stabilizer is infinite. Let \(\XX\) be 
  the Zariski closure of \(p\) in \(\GV\) and let \(a\) be a point realizing \(p\). Assume to the contrary, that 
  \(\Stab^{\mu}(p)\) is finite. Then \(\res(\mutube{p}\cdot{}a^{-1}\cap\GO)\) 
  is finite by Prop.~\ref{prp:Gp}. Note that \(\mu\cdot \XX\cdot a^{-1}\cap \GO=\bigcup_q (\mu \cdot q\cdot a^{-1}\cap \GO )\), where $q$ ranges over types centered at infinity on $\XX$. By Cor.~\ref{cor:fin}, there are only 
  finitely many types centered at infinity on \(\XX\), so the set 
  \(\XX\cdotp{}a^{-1}\cap\GO\) is the intersection of \(\XX\cdotp{}a^{-1}\) 
  with a (disjoint) union of finitely many balls \(\mu\cdotp{}g\), for 
  \(g\in\GR(\kk)\).
  
  Therefore, \(\XX\cdotp{}a^{-1}\cap\GO\) is a non-empty \(v{+}g\)-open 
  subset of \(\XX\cdotp{}a^{-1}\).  However, it is also a \(v{+}g\)-closed 
  subset of \(\XX\cdotp{}a^{-1}\) since \(\GO\) is \(v{+}g\)-closed. By 
  Prop.~\ref{prp:defcon}, \(\XX\cdotp{}a^{-1}\) is definably connected, so  
  \(\XX\cdotp{}a^{-1}\subseteq\GO\subseteq\AAA^n\).  However, this is 
  impossible since this implies that \(\XX\cdotp{}a^{-1}\) as an affine curve 
  is bounded in \(\AAA^n\), contradicting Theorem~\ref{thm:boundedvar}.
\end{proof}

\begin{rmk}[Relation to o-minimal PS-subgroups]\label{rmk:psgroups}
Recall that in the o-minimal context, for each definable group \(\GV\) and an 
unbounded semi-algebraic curve \(\gamma:(a,b)\ra{} \GV\), we use \(H_\gamma\) 
to denote the PS-subgroup of \(\gamma\) (or in other words, the o-minimal 
\(\mu\)-stabilizer of the type of \(\gamma\) at \(b\)).

In the case when \(\kk\) is \(\mathbb{C}\), the \(\mu\)-stabilizers of a point at infinity in Theorem~\ref{thm:main1a} is closely related to the group \(G_\gamma\), the Zariski-closure of \(H_\gamma\), as described in the construction in the 
beginning of the introduction.
(both viewed as definable in \(\mathcal{R}\))

Namely, assume we are given a complex curve \(\XX\) embedded in the complex 
affine algebraic group \(\GV\), and a point \(\alpha\in\ti{\XX}\setminus{}\XX\). 
Let \(\gamma:(a,b)\ra{}\XX\) be a semi-algebraic curve over \(\real\) whose 
limit at \(b\) is \(\alpha\) (in the sense discussed in the introduction).

Let \(\mathcal{R}\) be a sufficiently saturated real closed field extending 
\(\real\), and let \(\mathcal{C}\) be \(\mathcal{R}^2\), viewed as an 
algebraic closure of \(\mathcal{R}\). We may view \(\XX\), \(\GV\) and \(\gamma\) 
as definable in \(\mathcal{R}\). By~\cite{Ser}, one can compute the 
PS-subgroup \(H_\gamma\) associated to \(\gamma\) as the \(\mu\)-stabilizer of 
\(tp_{sa}(\alpha/\mathbb{R})\), the type of \(\alpha\) in the theory of real 
closed fields.

Let \(\OO_R\) be the convex hull of \(\real\) in \(\mathcal{R}\), and let 
\(\OO={\OO_R}^2\), viewed as a subring of \(\mathcal{C}\). Then 
\(\mathcal{C}\) equipped with \(\OO\) as a valuation ring is a model of 
\(T_{\mathrm{loc}}\), and the type \(q=tp_{\ACVF}(\alpha/\mathbb{C})\) is 
contained in \(p=tp_{sa}(\alpha/\mathbb{R})\) because the maximal ideal can 
be viewed as a partial type.  Hence \(H_\gamma\subseteq \Stab^\mu(q)\).  
However, \(\Stab^\mu(q)\) is 1-dimensional and \(H_\gamma\) is infinite, so 
\(G_\gamma \subseteq \Stab^\mu(q)\) and the index is finite.
\end{rmk}

\section{Proof of the main theorem}\label{sec:general}
\subsection{Residually algebraic saturation}
We would like to work with saturation in a residually algebraic context, 
i.e., without extending the residue field. Thus we make the following 
definition.

\begin{Def}
  A model \(K\) of \(T_{loc}\) is (sufficiently) \deft{\(\VG\)-saturated} 
  if every \(\LL_\val\)-residually algebraic type over a (sufficiently) small subset of 
  \(K\) is realised in \(K\).
\end{Def}
\begin{thm}
  Let \(L\) be a (possibly trivially) valued field, then there is a \(\VG\)-saturated extension of \(L\).
\end{thm}

\begin{proof}
Let \(\Gamma\) be a sufficiently saturated ordered abelian group and \(k\)  the algebraically closed closure of \(\RES(L)\). Consider the Hahn series field
  \[
    k((t^\Gamma))=\St{\sum_\gamma c_\gamma t^\gamma}
    {c_\gamma\in k, \{\gamma:c_\gamma\neq 0\}\text{ is well ordered }}
  \]
  Clearly \(L\) embeds into \(K\) (see~\cite{kaplansky} for example).
Then, by a result of Poonen 
  (\cite[Theorem 2]{poonen}), 
  is a \(\VG\)-saturated model (with residue field \(k\)).
\end{proof}

\emph{From now on, \(K\) will be a fixed sufficiently \(\VG\)-saturated model \(K\)
with residue field \(\kk\) and we will identify definable sets and \(p\in 
S_G(\kk)\) with their realizations in \(K\). In particular, we will only 
consider residually algebraic types unless otherwise stated.}

As a first application, we note:
\begin{lem}
  Let \(p\in{}S_G(\kk)\) be residually algebraic. Then there is
  \(q\in{}p_\mu\) which is \(\mu\)-reduced and residually algebraic.
\end{lem}
\begin{proof}
  Let \(a\) be a realization of \(p\) in \(K\). There is a variety \(\VV\) 
  over \(\kk\) of minimal dimension that intersects \(\mu\cdot{}a\). The 
  above can be expressed as \(\LL_\val\)-formula, so it is witnessed by some element of \(K\).  
  Take \(q\) to be the \(\LL_\val\)-type of this element over \(\kk\).
\end{proof}

We would like to give a syntactic (or geometric) description of types 
realised in \(K\). To this end, we use another result from~\cite{HL}, which 
requires the following definition.

The following is a part of Lemma~9.1.1 in~\cite{HL}, which will be needed in the proof.
\begin{lem}(Hrushovski, Loeser)\label{lem:gpair}
  Let \(F\) be a valued field, \(\VV\) be an \(F\)-variety, and \(\XX\subseteq \VV\) be a 
  \(F\)-definable \(g\)-open set, then \(\XX(M_2)\subseteq \XX(M_1)\) 
  whenever \(M_1\) and \(M_2\) are algebraically closed valued field extensions of \(F\) with the same underlying field, and \(\OO_{M_1}\subseteq \OO_{M_2}\).
\end{lem}

We now have the following description:
 
\begin{prop}\label{sat}
 Let \(\Phi(x)\) be a small finitely consistent collection of \(g\)-open 
 sets, with parameters in \(L\subseteq{}K\). Then \(\Phi\) is realised in 
  \(K\). In addition, if \(p\) is a \(\LL_\val\)-residually algebraic type, 
  then it is the intersection of the \(g\)-open formulas that it implies.
\end{prop}

In other words, every partial type \(\Sigma\) of \(g\)-open sets admits an extension to a
\(\LL_\val\)-residually algebraic complete type \(p\) over the same set of 
parameters.

\begin{rmk}
  It is worth pointing out that Proposition \ref{sat} has an easy proof in 
  the case  when \(L=\kk=\complex\). See \cite[Section 3.1]{PS-flow}.
\end{rmk}

\begin{proof}
  Let \(b\) be any realisation of \(\Phi\) in \(\Uu\), and let \(k\) be the 
  residue field of \(L(b)\). Then \(k\) is the function field of some variety 
  \(\XX\) over \(\kk\), fix a valuation \(\val'\) of \(k\) over \(\kk\), with 
  residue field \(\kk\).

  Let \(M_2\) be the algebraic closure of \(L(b)\) with the induced valuation 
  from \(\Uu\). Consider a valuation of \(\RES(M_2)\) extending \(\val'\).  
  Abusing notaion, we call the valuation \(\val'\) as well. Let \(\bar{\OO}\) 
  be the valuation ring of \(\RES(M_2)\). Consider 
  \(\res^{-1}(\bar{\OO})\subseteq M_2\), this is again a valuation ring of 
  the underlying field of \(M_2\) over \(\kk\). We use \(M_1\) to denote the 
  same field as \(M_2\) with the valuation determined by 
  \(\res^{-1}(\bar{\OO})\).  Note that \(M_1\) has residue field \(\kk\).

  Then by Lemma~\ref{lem:gpair}, \(\phi(M_2)\subseteq\phi(M_1)\) for each 
  \(\phi\in\Phi\). In particular, \(b\) is a realisation of \(\Phi\) in 
  \(M_1\). But the residue field of \(M_1\) is \(\kk\), so \(tp_{M_1}(b/L)\) 
  is residually algebraic and hence realizable in \(K\).

  For the converse, let \(p\) be a complete \(\LL_\val\)-residually algebraic 
  type. By quantifier elimination in \(\ACVF\), it is given within its 
  Zariski closure by formulas of the form \(f(x)\ne{}0\), 
  \(\val(f(x))>\val(g(x))\) and \(\val(f(x))=\val(g(x))\ne\infty\). Each 
  formula of the last form is equivalent to \(\val(f(x)/g(x))=0\), so that 
  \(f(x)/g(x)\) has non-zero residue. Since \(p\) is residually algebraic, 
  the residue is actually a well determined element \(b\) of \(\kk\), so the 
  original formula is implied by \(\val(f(x)-bg(x))>\val(bg(x))\), which is also in $p$.
\end{proof}

We now apply this result in our context:

\begin{cor}\label{cor:resalgtube}
  Let \(p\in S_G(\kk)\) and residually algebraic. Then 
  \(\mu(K)\cdotp{}p(K)=(\mutube{p})(K)\).
\end{cor}
\begin{proof}
 Since \(K\) is contained in the monster model, 
 \(\mu(K)\cdotp{}p(K)\subseteq(\mutube{p})(K)\). For the reverse containment, 
 for \(p\) residually algebraic, fix \(a\in(\mutube{p})(K)\).
 Recall that it means that for any \(\phi\in{}p\), there is
 \(\ve_\phi\in\mu\) such that \(\models\phi(\ve_\phi\cdotp a)\).  Since \(p\) 
 is residually algebraic, we may, by Prop.~\ref{sat}, assume that each such 
 \(\phi\) is \(g\)-open.

 Consider the following partial type: 
 \(
   \Sigma(y)=\St{\phi(y\cdotp{}a)\land\mu(y)}{\phi{\in}p\,\,\text{\(g\)-open}}
 \).
 Each \(\phi\) there is \(g\)-open, hence also \(\phi(y\cdotp{}a)\) (since 
 the group is algebraic) and \(\mu\) is given by strict inequalities, so this is a small collection of \(g\)-open sets, consistent by assumption. By the other  direction of Prop.~\ref{sat}, we can find \(\ve\in\mu(K)\) such that \(\ve\cdotp{}a\) satisfies \(p\).
\end{proof}

\begin{cor}\label{model}
 Let \(p\) be a residually algebraic \(G\)-type over \(\kk\), and let \(a\) 
 be a realization in \(K\). Then 
 \(\Stab^\mu(p)(\kk)=\res(\mu(K)\cdotp{}p(K)\cdotp{}a^{-1}\cap{}\GO)\).
\end{cor}
\begin{proof}
  Since \(\res(\mu(K)\cdotp{}p(K)\cdotp{}a^{-1}\cap{}\GO)\subseteq \res(\mu(\Uu)\cdotp{}p(\Uu)\cdotp{}a^{-1}\cap{}\GO)\), we have that \(\Stab^\mu(p)(\kk)\supseteq \mu(K)\cdotp{}p(K)\cdotp{}a^{-1}\cap{}\GO\) by Cor.~\ref{def}. The reverse containment follows from Cor.~\ref{cor:resalgtube}.
\end{proof}

\subsection{\texorpdfstring{\(\mu\)}{mu}-reduced types and their stabilizers}
In this section we prove Cor.~\ref{cor:fing}, an analogue of 
Corollary~\ref{cor:fin} for types of higher dimension.

\emph{Recall that we are working within \(K\), a \(\Gamma\)-saturated model, 
  and all the elements in the statement are from \(K\) and definable sets are 
  identified with their realization in \(K\)}.

  In particular, we have the following.
  \begin{lem}\label{lem:redinv}
  If \(a\) is \(\mu\)-reduced and \(g\in\GO\), then \(g\cdot a\) is also 
  \(\mu\)-reduced, of the same dimension.
\end{lem}
\begin{proof}
  Assume \(\ve{}\cdot g\cdot a\in\WW\) with \(\ve\in\mu\) and \(\WW\) a variety over 
  \(\kk\). Since \(\ve{}\cdot g\cdot a=\bar{g}\cdot \ve'\cdot a\) for some \(\bar{g}\in\GR(\kk)\) and 
  \(\ve'\in\mu\), we have \(\ve'\cdot a\in\bar{g}^{-1}\cdot\WW\), a variety over \(\kk\) 
  of the same dimension as \(\WW\).
\end{proof}
The following is an important observation about \(\mu\)-reduced types.

\begin{prop}\label{prp:comptype}
  Let \(p\in{}S_G(\kk)\) be a \(\mu\)-reduced residually 
  algebraic type centered at infinity, and let \(a\models{}p\). Let \(\VV\) be 
  the unique irreducible \(\kk\)-variety such that \(a\in\VV\) and 
  \(\dim(\VV)=\dim(p)\). For any definable set \(\XX\), assume that \(\XX\subseteq \GO\cdotp{}a\cap\VV\) be definably 
  connected and \(a\in \XX\). Then for every \(b\in\XX\) we have \(tp(a/\kk)=tp(b/\kk)\), where \(tp(\cdot/\kk)\) denotes the \(\LL_\val\)-type over \(\kk\).
\end{prop}

\begin{proof}[Proof of~\ref{prp:comptype}]
  By Lemma~\ref{lem:redinv}, \(b\) is not contained in any proper 
  subvariety of \(\VV\), so is nonzero when evaluated by any regular function on \(\VV\).  
  Hence, every element of the function field \(\kk(\VV)\) is well defined as 
  a \(\kk\)-definable function on \(\XX\).

  Assume that the types of \(a\) and \(b\) are different. By quantifier 
  elimination in \(\ACVF\), without loss of generality, there is \(f\in\kk(\VV)\) such that 
  \(\val(f(a))<0\le\val(f(b))\). We may further assume that the last 
  inequality is strict, by subtracting the residue.

  By \cite{HL}, it can be easily checked that rational functions are 
  \(v{+}g\)-continuous on their domain, so the image \(f(\XX)\) is again 
  definably connected. As a definable subset of \(K\), it is a union of 
  ``Swiss cheeses", and by definable connectedness, the Swiss cheese decomposition of the image will be of the form
  \(\BB\setminus\cup_{i\leq{}m}\CC_i\), where \(\BB\) is a ball and \(\CC_i\)'s are disjoint 
  sub-balls of \(\BB\).

  \begin{claim*}
    \(f(\XX)\) contains a \(\kk\)-point.
  \end{claim*}
  \begin{proof}[Proof of claim]
    Since \(\BB\) contains both a point with positive valuation and point with valuation \(\leq 0\), then it must contain \(\OO\). If \(f(\XX)\) contains no \(\kk\)-point, \(\kk\) must be covered by \(\bigcup_{i\leq{}m}\CC_i\). This implies one of the \(\CC_i\) contains at least two points in \(\kk\) and hence contains \(\OO\). But this is a contradiction since that means that there is no point in 
    \(f(\XX)\) with positive valuation.
  \end{proof}
  Hence, we know that there must be some \(c\in\kk\) such that 
  \(c\in{}f(\XX)\). Note however that each element in \(\XX\) is a generic 
  point of \(\VV\) by Lemma \ref{lem:redinv} and we know that this would imply that the rational 
  function \(f\) is constant, a contradiction to the assumption. Hence we 
  know that \(tp(a/\kk)=tp(b/\kk)\).
\end{proof}

Using a similar argument, we have the following, which is the key fact that 
will replace Cor.~\ref{cor:fin} for our proof of the main theorem.

\begin{cor}\label{cor:fing}
  Let \(a\in\VV(K)\) be \(\mu\)-reduced, with \(\VV\) the Zariski closure of 
  \(a\) over \(\kk\). Then there are finitely many types \(p_1,...,p_m\in S_G(\kk)\) for some \(m\) such that if \(g \in \GV(\OO)\) and 
  \(g\cdotp{}a\in\VV\), then \(tp(g\cdotp{}a/\kk)=p_i\) for some \(i\).
\end{cor}
\begin{proof}
  From Theorem~\ref{comp}, we know that there are only finitely definably 
  connected components of the set \(\GO\cdotp{}a\cap\VV\), call them 
  \(\XX_i\) for \(i=1,..., n\) for some \(n\). By Proposition \ref{prp:comptype}, we have that for each \(b,b'\in \XX_i\), we have that 
  \(tp(b/\kk)=tp(b'/\kk)\). Hence there are only finitely many types 
  \(p_i\)'s with the property stated in the statement of the corollary.
\end{proof}

Here, we stated a variant that is similar to the \cite{Ser}'s Claim~3.13.

\begin{cor}\label{sep}
  In the same setting as above, there is a \(\LL_\val\)-definable set \(\XX\) over \(\kk\)
  containing \(a\), such that for each \(b\in\GO\cdotp{}a\cap\XX\), 
  \(tp(b/\kk)=tp(a/\kk)\). Furthermore, \(\XX\) is \(v{+}g\)-open and 
  \(\XX\cdotp{}a^{-1}\) is \(v{+}g\)-clopen in \(\GO\).
\end{cor}
\begin{proof}
  By the above corollary, we see that there are finitely many regular functions \(f_{ij},g_{ij}\) such that for the formula \(\val(f_{ij})<\val(g_{ij})\) \(p_i\) and \(p_j\) disagrees. The set defined by some boolean combinations of the above formulas containing \(tp(a/\kk)\) will define the set \(\XX\).
\end{proof}

In particular, we have the following.
\begin{cor}\label{cor:main1}
\(\Stab^\mu(p)(\kk)\) is infinite for each \(p\) residually algebraic and 
centered at infinity.
\end{cor}
\begin{proof}
  Without loss of generality, we can assume that \(p\) is \(\mu\)-reduced and 
  let \(a\models{}p\) be any realization and \(\VV\) denote its Zariski 
  closure. We have \(\dim\VV>0\) since \(p\) is centered at infinity. Also, 
  \(\VV\) is an irreducible \(\kk\)-variety hence \(v{+}g\)-connected and so 
  is \(\VV\cdotp{}a^{-1}\). By Proposition~\ref{prp:Gp}, if 
  \(\Stab^\mu(p)(\kk)=\GR_p(\kk)\) is finite, \(\VV\cdotp{}a^{-1}\cap\GO\) 
  can be covered by finitely many \(v{+}g\)-open sets and hence is 
  \(v{+}g\)-open. But \(\VV\cdotp{}a^{-1}\setminus\GO\) is also 
  \(v{+}g\)-open by definition, a contradiction to the fact that non-zero 
  dimensional affine varieties are not bounded in the affine space.
\end{proof}

It is worth noting that the above proof uses the same idea in the 
\(1\)-dimensional case where the key ingredient is the connectedness of 
irreducible varieties.

\subsection{Dimension of the \texorpdfstring{\(\mu\)}{mu}-stabilizers}
\label{sec:dim}
Corollary~\ref{cor:main1} forms the first part of the main theorem. Before 
proving the other part, we need some machinery about varieties over \(\OO\).  
The main facts can be found in \cite{hrushovski2017valued} and 
\cite{halevi2016stably}.

\begin{Def}
  Let \(\OO\) be a valuation ring, let \(L=\mathrm{Frac}(\OO)\) and 
  \(k=\res(\OO)\). By a \deft{variety over \(\OO\)}, we mean a flat reduced 
  scheme \(\Vv\) of finite type over \(\OO\). In particular, it has a generic 
  fiber \(\VV_L\), which is a variety over \(L\), obtained by base change 
  with respect to the morphism \(\OO\to L\), and a special fiber \(\VV_k\), 
  which is a variety over \(k\), obtained by base change with respect to 
  \(\OO\to k\).
\end{Def}

\begin{rmk}
  It is worth noting that since \(\OO\) is a valuation ring, 
  \(A=\OO[x_1,...,x_n]/I\) is flat over \(\OO\) iff no nonzero element in 
  \(\OO\) is a zero divisor in \(A\).

  In particular, if \(S\) is any subset of \(\OO^n\), then 
  \[
  I=\St{f\in\OO[x_1,...,x_n]}{f(s)=0\,\forall{}s\in{}S}\subseteq\OO[\bar{x}]
  \]
  is an ideal and \(A=\OO[x_1,...,x_n]/I\) is flat over \(\OO\).

  We use \(I_L\) and \(I_k\) to denote the ideal generated by \(I\) in
  \(A\otimes L\) and \(A\otimes k\) respectively. Then the generic 
  (respectively special) fiber of \(\mathrm{Spec}(A)\) is 
  \(\mathrm{Spec}(A\otimes L/I_L)\) (respectively \(\mathrm{Spec}(A\otimes 
  k/I_k))\). Given an affine variety \(\VV\) over \(L\), we may always choose 
  a variety \(\Vv\) affine over \(\OO\) whose generic fibre is \(\VV\).
\end{rmk}

The following is~\cite[Theorem~3.2.4]{halevi2016stably}.

\begin{thm}[Halevi]\label{thm:halevi}
  Let \(K\) be a model of ACVF and \(\Vv\) be an irreducible variety over 
  \(\OO_K\). If \(\VV_K\) has an \(\OO_K\)-point then the \(\OO_K\)-points 
  are Zariski dense, and the canonical map
  \(\res:\VV_K(\OO_K)\rightarrow\VV_k(k)\) is surjective, where \(\res\) is 
  given by taking residue pointwise.
\end{thm}

Now let us get back to the proof of the main theorem, we will prove the part 
concerning the dimension now.

\begin{proof}[Proof of Theorem~\ref{thm:main}]
  The first part of the theorem is the statement of 
  Corollary~\ref{cor:main1}, so it remains to show the equality of 
  dimensions, in the case the type \(p\) is \(\mu\)-reduced. Under this 
  assumption, let \(a\) realize \(p\), and let \(\VV\) be its Zariski 
  closure. In particular, \(\VV\) is irreducible.

  We first note that if \(g\in\GV(\OO)\) and \(g\cdot a\in\VV\), then 
  \(tp(g\cdot a/\kk)\) will be one of the finitely many types \(p_i\) 
  provided by Corollary~\ref{cor:fing}.
  Hence \(\mu\cdot \VV\cdotp{}a^{-1}\cap\GO\) will be a finite union of 
  cosets of \(\mu\cdot \GV_p\). Thus it suffices to show that 
  \(\res(\VV\cdotp{}a^{-1}\cap\GO)\) and \(\VV\) have the same dimension.

  As above, the affine \(K\)-variety \(\VV\cdotp{}a^{-1}\) can be viewed as 
  the generic fibre of a variety \(\Vv\) over \(\OO_K\).  Furthermore, 
  \(\Vv\) has an \(\OO_K\)-point, namely \(e\), the identity of the group 
  \(G\). It follows from Theorem~\ref{thm:halevi} that the map 
  \(\res:\VV\cdotp{}a^{-1}\cap\GO\rightarrow\GR\) is onto the special fiber 
  of \(\VV\cdotp{}a^{-1}\).  Also, by flatness, the special fiber has the 
  same dimension as the generic fiber, which is the dimension of \(p\).

  This completes the proof of the dimension part in the main result. The 
  solvability of the group is proved separately, as Theorem~\ref{thm:main3}.
\end{proof}

The proof also shows that the special fiber, being the image of 
\(\VV\cdotp{}a^{-1}\cap\GO\), is a finite union of cosets of 
\(\Stab^\mu(p)\). Therefore, we have established the following.

\begin{cor}
  Let \(\VV\subseteq \GV\) be a variety over \(\kk\), let \(a\models{}p\), 
  where \(p\in S_G(\kk)\) is a \(\mu\)-reduced residually algebraic type 
  centered at infinity. Assume further that the Zariski closure of \(p\) over 
  \(\kk\) is \(\VV\).
  
  Then the special fiber of \(\VV\cdotp{}a^{-1}\) is equi-dimensional i.e.  
  each irreducible component of it has the same dimension. Moreover, each 
  irreducible component of the special fiber of \(\VV\cdotp{}a^{-1}\) is a 
  coset of an algebraic subgroup of \(\GR\).
\end{cor}

\subsection{Structure of \texorpdfstring{\(\Stab^\mu(p)\)}{Stab}}
In this section, we analyze the structure of \(\Stab^\mu(p)\). Note that due 
to trivial constraints on characteristic, it is not possible to show in 
general such a group is torsion-free. However, in characteristic \(0\), we 
can indeed show it is torsion free.

\begin{lem}\label{sol}
  Let \(p\in{}S_G(\kk)\) be residually algebraic and let \(H\) be a 
  \(\kk\)-definable linear subgroup of \(G\) with \(p\in{}H\). Then 
  \(\Stab^\mu(p)\) computed in \(\GV\) and in \(\HH\) coincide, where \(\HH\) denotes the group \(H\) viewed as a subset in \(\VF\).
\end{lem}
\begin{proof}
  Since the Zariski closure \(\VV\) of \(p\) is contained in \(\HH\) in this case since \(\HH\) is a Zariski-closed subgroup and \(\mu_G\cap \HH=\mu_H\). Hence, the arguments of computing the \(\mu\)-stabilizers of \(p\) can be carried out in both \(\HH\) and \(\GV\) and the results will be the same.
\end{proof}

The following is the Iwasawa Decomposition over non-archimedean fields, it 
can be found in~\cite[Proposition 4.5.2]{bump}.
\begin{thm}
  Let \(\GV\) be a reductive linear algebraic group over \(\kk\), there is a 
  solvable subgroup \(H\) over \(\kk\) such that 
  \(\GV(K)=\GV(\OO)\cdotp{}\HH(K)\).
\end{thm}

For \(\GL_n\) we may take \(H\) to be the standard Borel subgroup (upper triangular 
matrices).

\begin{thm}\label{thm:main3}
  Let \(p \in S_G(\kk)\) be centered at infinity and residually 
  algebraic. Then \(\Stab^\mu(p)\) is solvable.
\end{thm}
\begin{proof}
  We can embed \(\GV\subseteq\GL_n\) over \(\kk\) for some \(n\), and use the 
  Lemma~\ref{sol} to reduce to the case \(\GV=\GL_n\). Let \(\HH\) be the 
  Borel.  By the Iwasawa decomposition, we have some \(g\in\GV(\OO)\) such 
  that \(g^{-1}\cdotp{}a=\beta\in{}\HH(K)\). Let \(g_1\in\GR(\kk)\) be such 
  that \(g_1\cdotp g^{-1}\in\mu\). Hence 
  \(g_{1}^{-1}\cdotp{}a\in\mu(K)\cdotp\beta\), so 
  \(\Stab^\mu(g_1^{-1}\cdotp{}p)=\Stab^\mu(q)\subseteq{}\ol{\HH}\). By 
  Lemma~\ref{lem:conj}, this group is conjugate to \(\Stab^\mu(p)\), hence 
  \(\Stab^\mu(p)\) is solvable.
 \end{proof}

\begin{cor}
If \(G\) is not solvable and \(G\) is irreducible, then there is no \(\mu\)-reduced residually algebraic \(G\)-type of full dimension. 
\end{cor}

\begin{rmk*}
 We briefly introduce the Zariski-Riemann space of a variety over \(k\), and explain its connection with our setting.
\begin{Def*}
Let \(V\) be an variety over \(k\), the Zariski-Riemann space of \(V\) over \(k\), is the set of valuation rings of \(k(V)\) over \(k\), denoted by \(\mathbf{RZ}_k(V)\).
\end{Def*}
Note that by quantifier elimination in ACVF, for a linear algebraic group \(\GV\) over \(\kk\), it is not hard to see that the above set \(\mathbf{RZ}_\kk(\GV)\) is in canonically embeddable into the set \(S_G(\kk)\). Hence we can identify \(\mathbf{RZ}_\kk(\GV)\) with its image in \(S_G(\kk)\).
Note further, since \(\mu\) is Zariski dense in \(\GV\), we see that for each 
\(p \in S_G(\kk)\), there is some \(q \in \mathbf{RZ}_\kk(\GV)\) such that 
\(p\sim_\mu q\).

The above argument implies that we can consider the quotient of 
\(\mathbf{RZ}_\kk(\GV)\) under \(\mu\), which exactly the space 
\(S^\mu_G(\kk)\). Note further that the equivalence relation induced by 
\(\mu\) on \(\mathbf{RZ}_\kk(\GV)\) is independent of the 
\(\kk\)-closed-immersion of \(\GV\) into \(\AAA^n\), since every embedding 
over \(\kk\) will respect \(\mu\). We will explore more on the relation 
induced by \(\mu\) on \(\mathbf{RZ}_\kk(\GV)\) in a subsequent paper.
\end{rmk*}
\bibliographystyle{plain}
\bibliography{bibliography}

\end{document}